\documentclass[10pt,a4paper]{article}

\usepackage[pdftex]{graphicx}
\usepackage{multirow}
\usepackage{morefloats}     
\usepackage{url}    
\usepackage{algorithm}
\usepackage{amsthm}         
\usepackage{amsmath}
\usepackage{colortbl}   

\include{Definitions}     

\newcommand{\floor}[1]{\lfloor #1 \rfloor}
\newcommand{\ceil}[1]{\lceil #1 \rceil}

\theoremstyle{definition}
\newtheorem{defn}{Definition}[section]

\begin{document}

\title{Prime Sums of Primes}

\author{Dmitry Kamenetsky \\
dkamenetsky@gmail.com \\
Adelaide, Australia
}

\maketitle

\abstract{We present a variety of prime-generating constructions that are based on sums of primes. 
The constructions come in all shapes and sizes, varying in the number of dimensions and number of generated primes.
Our best result is a construction that produces 6 new primes for every starting prime.
}

\section{Introduction}

Constructions made from primes have fascinated mathematicians for many decades due to the beauty of their design.
A number of such constructions have been proposed, such as: prime magic squares~\cite{Magic,Magic2},
prime arrays~\cite{Array} and primes in arithmetic progressions~\cite{Progression,Progression2}.

In this paper we investigate some new prime-generating constructions that are based on sums of primes.
Our constructions come in two flavours: \emph{standard} and \emph{recursive}. In standard constructions
new primes are generated as the sum of primes used in the construction.
Recursive constructions generate new primes, which
in turn generate further primes. The recursion terminates when no more primes can be generated.
Typically we only use odd primes (ignore 2), forcing our sums to contain an odd number of elements.
Our overall aim is to generate constructions of the largest size (\emph{order}).
If two constructions have the same order then we typically prefer the one with smallest sum of elements (\emph{weight}).
To find all the constructions we use a variant of the randomised hill-climbing algorithm.
For small constructions we were able to find the optimal solutions (smallest weight) by using a brute force method.

We describe the following standard constructions: prime vectors (Section~\ref{sec:vectors}), cyclic
prime vectors (Section~\ref{sec:cyclicVectors}), Goldbach squares (Section~\ref{sec:goldbach}) and 
prime matrices (Section~\ref{sec:matrices}). We describe the following recursive constructions:
prime tuples (Section~\ref{sec:tuples}), prime stairs (Section~\ref{sec:stairs}),
prime pyramids (Section~\ref{sec:pyramids}) and prime cylinders (Section~\ref{sec:cylinders}).

\section{Prime Vectors}
\label{sec:vectors}

\begin{defn}
A \emph{prime vector} of order $n$ is an array of distinct primes $P=(p_0,p_1, \ldots, p_{n-1})$, such that every sum of an odd number of 
consecutive elements is also prime. In other words
\begin{equation}
\begin{split}
\sum_{0 \leq k \leq 2L} P(i+k) & ~~~\mbox{is prime for} \\
& \forall~i~\mbox{such that}~0 \leq i \leq i+2L < n.
\end{split}
\end{equation}
\end{defn}

In the above definition, $i$ is the index of the first prime in each sum, while $(2L+1)$
is the number of terms in each sum.\footnote{If $L=0$ then we have a singleton rather than a sum.}
For a given $n$ there are $\floor{(n-1)^2/4}$ sums. 
Consider a prime vector of
order 5: $(3,11,5,7,17)$. Its every element is prime, as well as, every sum of an odd number of consecutive elements:

\begin{equation}
\begin{split}
3+11+5=19,~~~11+5+7=23, \\
5+7+17=29, ~~~3+11+5+7+17=43.
\end{split}
\end{equation}

We used a variant of hill-climbing to find prime vectors (see Algorithm~\ref{alg:array}).
We start with a random array of distinct primes and then perform various mutations, such as swapping two primes or replacing one prime with
a new one. If the mutation improves the score then we keep it, otherwise we revert it. The score measures the number of ``incorrect" (composite)
sums that the array generates. Hence we want to minimise this score. Using this algorithm we were able to obtain a prime vector of order
23 that generates 121 primes\footnote{Prime vectors of smaller orders are sub-arrays of this array.}: \\

(239, 131, 109, 181, 83, 43, 41, 223, 53, 233, 271, 103, 269, 71, 19, 47, 241, 23, 277, 199, 281, 29, 37).\\

For small orders it is possible to obtain multiple solutions. In such cases we choose the solution with the smallest \emph{weight} - sum of
all elements. In fact, this allows us to define an \emph{optimal prime vector}:
\begin{defn}
A prime vector is \emph{optimal} if its weight is the lowest possible.
\end{defn}

\begin{algorithm}[!htpb]
\renewcommand{\arraystretch}{1.15}
\caption{: Algorithm for finding prime vectors.}
\vspace{0.5ex}
\begin{tabular}{rl}
\\
1 & $bestScore \leftarrow \infty$ \\
2 & $S^* \leftarrow$ random set of n distinct primes \\
3 & \\
4 & while True\\
5 & ~~$S \leftarrow S^*$\\
6 & ~~mutate($S$) \\
7 & ~~$score \leftarrow$~score($S$) \\
8 & ~~if $score < bestScore$ \\
9 & ~~~~$bestScore \leftarrow score$ \\
10 & ~~~~$S^* \leftarrow S$ \\
11 & ~~~~print($S^*$) \\
12 & ~~end \\
13 & end \\
\end{tabular}
\label{alg:array}
\end{algorithm}

For $n \leq 14$ we were able to find the optimal prime vectors (see Table~\ref{tab:array1}). To achieve this we used a brute force algorithm.
This algorithm iterates through every permutation of $n$ distinct odd primes whose weight is below the best known weight. If a permutation forms
a prime vector then the best known weight is updated and the array is printed out. The algorithm terminates when there are no more permutations
whose weight is less than the best known weight. Table~\ref{tab:array1} also shows the running time of this algorithm.

For $n > 14$ we used Algorithm~\ref{alg:array} to find the upper bounds on the minimal weight (see Table~\ref{tab:array2}).
To obtain the lower bound we used sequences from the OEIS~\cite{OEIS}.
For odd $n$ the weight must be a prime, so we used sequence A068873 - smallest prime which is a sum of $n$ distinct primes.
For even $n$ we used sequence A071148 - sum of the first $n$ odd primes.


\begin{table}[!htpb]
\centering
\begin{tabular}{|c|c|c|c|}
\hline
$n$ & Prime Vector & Weight & Time \\ \hline
1   & $(2)$   & 2      &      \\ \hline     
2   & $(3,5)$ & 8      &      \\ \hline     
3   & $(3,5,11)$ & 19 &      \\ \hline      
4   & $(3,5,11,7)$ & 26 &      \\ \hline
5   & $(3,11,5,7,17)$ & 43 &      \\ \hline
6   & $(3,11,5,7,17,13)$ & 56 &      \\ \hline
7   & $(3,17,23,7,11,13,5)$ & 79 &      \\ \hline
8   & $(3,11,17,13,29,19,5,7)$ & 104 &      \\ \hline
9   & $(7,17,13,23,11,3,29,5,19)$ & 127 &      \\ \hline
10  & $(3,7,19,11,13,23,31,5,37,17)$ & 166 &      \\ \hline
11  & $(3,23,41,19,11,13,17,7,5,31,53)$ & 223 & 17s     \\ \hline
12  & $(7,41,19,11,23,3,5,29,13,47,43,17)$ & 258 & 8m      \\ \hline
13  & $(13,53,7,23,11,3,29,5,19,17,43,47,37)$ & 307 & 73m      \\ \hline
14  & $(17,43,47,13,29,5,3,23,11,19,41,7,53,37)$ & 348 & 14h      \\ \hline
\end{tabular}
\caption{Optimal prime vectors for $n \leq 14$, their weight and the time required to compute them. Computation times less than 1 second are
not shown.}
\label{tab:array1}
\end{table}

\begin{table}[!htpb]
\centering
\begin{tabular}{|c|c|c|c|c|c|c|c|c|c|}
\hline
$n$         & 15  & 16  & 17  & 18  & 19  & 20   & 21   & 22   & 23 \\ \hline
lower bound & 379 & 438 & 499 & 566 & 643 & 710  & 809  & 872  & 983 \\ \hline
upper bound & 443 & 522 & 641 & 888 & 983 & 1430 & 1627 & 1824 & 3203  \\ \hline
\end{tabular}
\caption{Best bounds on the minimal weight of prime vectors for $15 \leq n \leq 23$.}
\label{tab:array2}
\end{table}

\subsection{Cyclic Prime Vectors}
\label{sec:cyclicVectors}

We can also introduce a \emph{cyclic prime vector} and define its optimality in a similar fashion:
\begin{defn}
A \emph{cyclic prime vector} of order $n$ is a prime vector $P$ of order $n$ with the additional property that prime sums can span from the end to the
start of the array. In other words
\begin{equation}
\begin{split}
\sum_{0 \leq k \leq 2L} P((i+k)\bmod n) & ~~~\mbox{is prime for} \\
& \forall~i~\mbox{such that}~0 \leq i < n~\mbox{and} \\
& \forall~L~\mbox{such that}~0 \leq 2L < n.
\end{split}
\end{equation}
\end{defn}
For a given $n$ there are $(n-2)(2n-1+(-1)^n)/4$ sums.  
For example the cyclic prime vector $(5,7,17,13,11)$ generates the following 6 sums:

\begin{equation}
\begin{split}
5+7+17=29,~~~7+17+13=37, \\
17+13+11=41,~~~13+11+5=29, \\
11+5+7=23,~~~5+7+17+13+11=53.
\end{split}
\end{equation}

Cyclic prime vectors differ from normal prime vectors in a few key ways. 
Every cyclic prime vector is also a normal prime vector, but the opposite may not be the case.
Unlike normal prime vectors, cyclic prime vectors can be permuted without affecting their prime sums.
Also we cannot easily generate cyclic prime vectors as sub-arrays of larger cyclic prime vectors. Due to the cyclic
requirement, cyclic prime vectors require more prime sums for the same order, making them significantly harder to find.

Using the brute force algorithm described above we were able to find the optimal cyclic prime vectors for
$n \leq 10$ (see Table~\ref{tab:cyclicArray1}). The computation for the optimal
cyclic prime vector of order 11 was still running after 4 days, so it is not shown.
It is interesting to note that the weight for $n=9$ is smaller than the weight for $n=8$. Using an algorithm similar to
Algorithm~\ref{alg:array} we found cyclic prime vectors up to order 14 (see Table~\ref{tab:cyclicArray2}). The largest array
generates 84 primes.


\begin{table}[!htpb]
\centering
\begin{tabular}{|c|c|c|c|}
\hline
$n$ & Cyclic Prime Vector & Weight & Time \\ \hline
1   & $(2)$   & 2      &      \\ \hline
2   & $(3,5)$ & 8      &      \\ \hline
3   & $(3,5,11)$ & 19 &      \\ \hline
4   & $(5,7,17,19)$ & 48 &      \\ \hline
5   & $(5,7,17,13,11)$ & 53 &      \\ \hline
6   & $(5,29,7,11,19,37)$ & 108 &      \\ \hline
7   & $(5,7,17,13,29,31,11)$ & 113 &      \\ \hline
8   & $(11,17,43,47,13,19,29,31)$ & 210 &      \\ \hline
9   & $(7,17,13,11,19,41,29,37,23)$ & 197 & 9s     \\ \hline
10  & $(11,19,23,47,31,53,43,67,89,127)$ & 510 & 2m     \\ \hline
\end{tabular}
\caption{Optimal cyclic prime vectors for $n \leq 10$, their weight and the time required to compute them. Computation times less than 1 second are
not shown.}
\label{tab:cyclicArray1}
\end{table}


\begin{table}[!htpb]
\centering
\begin{tabular}{|c|c|c|}
\hline
$n$ & Cyclic Prime Vector & Weight \\ \hline
11   & $(23,73,17,13,71,19,11,193,59,137,67)$   & 683         \\ \hline
12   & $(73,47,43,137,97,59,151,239,31,163,89,131)$ & 1260    \\ \hline
13   & $(41,43,23,73,71,53,13,173,151,59,127,263,283)$ & 1373    \\ \hline
14   & $(73,179,41,97,43,197,199,173,379,311,131,37,29,421)$ & 2310   \\ \hline   
\end{tabular}
\caption{Smallest (by weight) cyclic prime vectors found for $11 \leq n \leq 14$.}
\label{tab:cyclicArray2}
\end{table}

\section{Prime Tuples}
\label{sec:tuples}

\begin{defn}
A \emph{prime tuple} of order $n$ (odd) with \emph{length} $k$ is an array of distinct odd primes $(p_0,p_1, \ldots, p_{k-1})$,
such that every term after the $n$-th term is the sum of the previous $n$ terms. In other words
\begin{equation}
p_i = \sum_{q=i-n}^{i-1} p_q, ~~~ \forall i \geq n.
\end{equation}
\end{defn}
Note it is sufficient to use the first $n$ terms to represent a prime tuple, since the remaining terms can be generated via
sums of previous terms. We seek to find prime tuples of order $n$ such that their length is greatest. For example, here is a
prime tuple of order $7$ with length 25 - the longest we have found:\\

(\textbf{157, 379, 487, 109, 13, 7, 271}, 1423, 2689, 4999, 9511, 18913, 37813, 75619, 150967, 300511, 598333, 1191667, 2373823, 4728733, 9419653, 18763687, 37376407, 74452303, 148306273).\\

The first $7$ terms are shown in bold. The \emph{weight} of a prime tuple of order $n$ is the sum of its first $n$ terms.
When two tuples of the same order have the same length, then we prefer the one with the smaller weight.

Table~\ref{tab:tuple} shows the best prime tuples that we found for $n \leq 19$. We have used a brute force approach to prove
that the prime tuples for $n \in \{3,5,9,11\}$ are optimal. We notice that for $n \mod 6=3$ and $n \mod 6=5$ the optimal
prime tuples have length $2n+1$ and must contain a 3. 

%

\begin{table}[!htpb]
\centering
\begin{tabular}{|c|c|c|c|}
\hline
$n$ & Prime Tuple & Length & Weight \\ \hline
3   & $(3,13,7)$   & 7      & 23     \\ \hline
\multirow{2}{*}{5} & $(17,3,19,7,13)$ & \multirow{2}{*}{11} & \multirow{2}{*}{59} \\
    & $(17, 5, 11, 23, 3)$ & & \\ \hline
7   & $(157, 379, 487, 109, 13, 7, 271)$ & 25   &  1423    \\ \hline    
9   & $(11, 47, 17, 23, 41, 5, 3, 13, 19)$ & 19 & 179     \\ \hline
11  & $(43, 7, 19, 13, 3, 17, 11, 5, 29, 41, 23)$ & 23 & 211  \\ \hline
\multirow{2}{*}{13}  & $(53, 137, 11, 17, 41, 227, 47,$ & \multirow{2}{*}{34} & \multirow{2}{*}{1163} \\
    & $101, 83, 5, 149, 263, 29)$ &  &     \\ \hline
\multirow{2}{*}{15}  & $(29, 5, 23, 11, 41, 47, 89, 17,$ & \multirow{2}{*}{31} & \multirow{2}{*}{491} \\
    & $71, 3, 7, 13, 37, 19, 79)$ &  &  \\ \hline
\multirow{2}{*}{17}  & $(5, 47, 53, 11, 17, 41, 89, 3, 61,$ & \multirow{2}{*}{35} & \multirow{2}{*}{647} \\
    & $43, 97, 19, 13, 7, 37, 31, 73)$ & &  \\ \hline
\multirow{2}{*}{19}  & $(89, 227, 29, 17, 5, 251, 269, 107, 101, 197,$ & \multirow{2}{*}{43} & \multirow{2}{*}{2081} \\
    & $41, 191, 173, 179, 47, 53, 71, 11, 23)$ & &  \\ \hline
\end{tabular}
\caption{Best prime tuples found for $n \leq 19$.}
\label{tab:tuple}
\end{table}

\section{Prime Stairs}
\label{sec:stairs}

\begin{defn}
A \emph{prime stair} of order $n \geq 3$ is a $\ceil{\frac{n}{2}} \times n$ matrix $P$ such that every element
$P(r,c)$ at row $r > 0$ and column $c$
is a distinct prime and each new row is generated from the previous row as follows:
\begin{equation}
P(r,c):=P(r-1,c-1)+P(r-1,c)+P(r-1,c+1).
\end{equation}
\end{defn}
For a given $r > 0$ we must have $c \in [r,n-r-1]$.
For a given $n$ there are $\floor{(n-1)^2/4}$ sums. 
As a shorthand we can represent a prime stair of order $n$ via its first (top) row only, i.e., using an array of
length $n$. For example, the prime stair $(13,17,7,5,11,3,23)$ looks like this:

\begin{table}[!htpb]
\centering
\begin{tabular}{c}
$13,17,7,5,11,3,23$ \\
$37,29,23,19,37$ \\
$89,71,79$ \\
$239$
\end{tabular}
\end{table}

The \emph{weight} of a prime stair is defined as the sum of all elements in the first row.
We were able to find the optimal prime stairs for
$n \leq 11$ (see Table~\ref{tab:stair1}). The computation for the optimal
prime stair of order 12 was still running after 4 days, so it is not shown.
Using an algorithm similar to
Algorithm~\ref{alg:array} we found prime stairs up to order 15 (see Table~\ref{tab:stair2}).
The largest stair generates 49 primes.

\begin{table}[!htpb]
\centering
\begin{tabular}{|c|c|c|c|}
\hline
$n$ & Prime Stair & Weight & Time \\ \hline
3   & $(3,5,11)$   & 19      &      \\ \hline
4   & $(7,5,11,3)$ & 26      &      \\ \hline
5   & $(7,13,11,5,3)$ & 39   &      \\ \hline
6   & $(7,17,13,11,5,3)$ & 56 &      \\ \hline
7   & $(13,17,7,5,11,3,23)$ & 79 &      \\ \hline
8   & $(5,17,31,11,19,7,3,13)$ & 106 &      \\ \hline
9   & $(29,23,37,13,11,5,7,19,17)$ & 161 & 8s      \\ \hline
10  & $(5,29,7,17,13,11,37,23,19,41)$ & 202 & 3m     \\ \hline
11  & $(7,17,19,5,73,11,13,23,31,29,41)$ & 269 & 2.5h     \\ \hline
\end{tabular}
\caption{Optimal prime stairs for $n \leq 11$, their weight and the time required to compute them. Computation times less than 1 second are
not shown.}
\label{tab:stair1}
\end{table}

\begin{table}[!htpb]
\centering
\begin{tabular}{|c|c|c|}
\hline
$n$ & Prime Stair & Weight \\ \hline
12   & $(37,13,17,23,7,67,5,59,19,61,29,11)$   & 348         \\ \hline
13   & $(29,19,11,127,89,7,17,37,5,31,23,43,41)$ & 479    \\ \hline
14   & $(53,17,67,29,13,5,19,149,31,101,79,11,7,43)$ & 624    \\ \hline
15   & $(433,139,491,97,89,163,29,7,5,61,17,79,263,541,83)$ & 2497   \\ \hline    
\end{tabular}
\caption{Best (by weight) prime stairs found for $12 \leq n \leq 15$.}
\label{tab:stair2}
\end{table}

\subsection{Prime Pyramids}
\label{sec:pyramids}

Similarly we can define a 3D version of the prime stair that we will call a \emph{prime pyramid}:

\begin{defn}
A \emph{prime pyramid} of order $n \geq 3$ is a $\ceil{\frac{n}{2}} \times n \times n$ matrix $P$ such that every element $P(k,r,c)$
at level $k > 0$, row $r$ and column $c$ 
is a distinct prime and each new level is generated from the previous level as follows:
\begin{equation}
P(k,r,c):=\sum_{-1 \leq dr \leq 1}~~\sum_{-1 \leq dc \leq 1} P(k-1,r+dr,c+dc).
\end{equation}
\end{defn}

For a given $k > 0$ we must have $r,c \in [k,n-k-1]$.
For a given $n$ there are $n(n-1)(n-2)/6$ sums. 
As a shorthand we can represent a prime pyramid of order $n$ via its first (bottom) level only, i.e., using a $n \times n$ array.
For example, Table~\ref{tab:pyramid} shows a prime pyramid of order 5:

\begin{table}[!htpb]
\centering
\begin{tabular}{ccc}
  Level 0 & Level 1 & Level 2 \\ \\
  \begin{tabular}{|c|c|c|c|c|}
  \hline
  73  & 11 & 67 & 71 & 53 \\ \hline
  101 & 41 & 43 & 79 & 83 \\ \hline 
  13  & 3  & 31 & 7  & 23 \\ \hline 
  17  & 61 & 37 & 5  & 29 \\ \hline 
  97  & 89 & 19 & 59 & 47 \\ \hline 
  \end{tabular}
  &
  \begin{tabular}{|c|c|c|}
  \hline
  383 & 353 & 457 \\ \hline 
  347 & 307 & 337 \\ \hline 
  367 & 311 & 257 \\ \hline 
  \end{tabular}
  &
  \begin{tabular}{|c|}
  \hline
  3119 \\ \hline
  \end{tabular}
\end{tabular}
\caption{Prime pyramid of order 5.}
\label{tab:pyramid}
\end{table}

The \emph{weight} of a prime pyramid is the sum of all elements in its first level.
We were able to find all the optimal prime pyramids up to order 8 (see Table~\ref{tab:pyramid2}).
We also found an order 9 prime pyramid with a weight of 27325, but its optimality is not confirmed (see Table~\ref{tab:pyramid9}).

\minrowclearance 0.25cm
\begin{table}[!htpb]
\centering
\begin{tabular}{|c|c|c|}
\hline
n & Prime Pyramid & Weight \\ \hline

3 &
\minrowclearance 0cm
\begin{tabular}{|c|c|c|}
\hline
7 & 11 & 29 \\ \hline 
19 & 17 & 23 \\ \hline 
5 & 13 & 3 \\ \hline 
\end{tabular}
& 127  \\[5ex] \hline

4 &
\minrowclearance 0cm
\begin{tabular}{|c|c|c|c|}
\hline
53 & 19 & 47 & 7  \\ \hline
37 & 3 & 41 & 13  \\ \hline 
43 & 5 & 29 & 17   \\ \hline
11 & 79 & 23 & 31   \\ \hline
\end{tabular}
& 458    \\[6ex] \hline

5 &
\minrowclearance 0cm
\begin{tabular}{|c|c|c|c|c|}
\hline
73 & 11 & 67 & 71 & 53   \\ \hline
101 &  41 &  43 & 79 & 83   \\ \hline
13 &  3 &  31 &  7 &  23   \\ \hline
17 & 61 &  37 &  5  & 29   \\ \hline
97 &  89  & 19  & 59 &  47   \\ \hline
\end{tabular}
& 1159  \\[8ex] \hline

6 &
\minrowclearance 0cm
\begin{tabular}{|c|c|c|c|c|c|}
\hline
97 & 43  & 47  & 149  & 79 &  3   \\ \hline
113  & 103 &  61 &  151 &  11  & 37   \\ \hline
109  & 53  & 107  & 19  & 127  & 67  \\ \hline 
31 &  23 &  101  & 29  & 13  & 7   \\ \hline
83  & 5 &  59  & 71  & 73  & 157   \\ \hline
137  & 139 &  41  & 89  & 131  & 17   \\ \hline
\end{tabular}
& 2582  \\[8ex] \hline

7 &
\minrowclearance 0cm
\begin{tabular}{|c|c|c|c|c|c|c|}
\hline
199  & 47  & 223  & 79  & 61  & 107  & 157   \\ \hline
229 &  89 &  5  & 71  & 29  & 163  & 211   \\ \hline
167  & 109  & 83  & 3  & 23  & 137  & 151   \\ \hline
97  & 197  & 43  & 73 &  59 &  19  & 31   \\ \hline
191  & 103  & 139  & 179 &  41 &  7  & 11   \\ \hline
17  & 173  & 227  & 37 &  13  & 149 &  127   \\ \hline
181  & 53  & 67  & 113  & 131  & 193  & 101  \\ \hline
\end{tabular} 
& 5115  \\[10ex] \hline

8 &
\minrowclearance 0cm
\begin{tabular}{|c|c|c|c|c|c|c|c|}
\hline
263  & 229  & 23  & 167  & 89  & 61  & 109  & 79   \\ \hline
97  & 149  & 173  & 67  & 281  & 211  & 59  & 47   \\ \hline
113  & 163  & 283  & 197  & 11 &  7  & 233  & 227   \\ \hline
277  & 19  & 293  & 223  & 181  & 251  & 43  & 241   \\ \hline
179  & 139  & 191 &  239  & 193 &  71  & 41  & 103   \\ \hline
101  & 107  & 13  & 83 &  73  & 137  & 269  & 37   \\ \hline
3  & 127  & 311  & 271  & 157  & 307  & 199  & 53   \\ \hline
5  & 17  & 313  & 131  & 257  & 151 &  29 &  31   \\ \hline
\end{tabular}
& 9204  \\[12ex] \hline
\end{tabular}
\caption{Optimal prime pyramids for $3 \leq n \leq 8$.}
\label{tab:pyramid2}
\end{table}

\clearpage

\begin{table}[!htpb]
\centering
\begin{tabular}{|c|c|c|c|c|c|c|c|c|}
\hline
11 &  79  & 349 &  461 &  433 &  859 &  683 &  587  & 631   \\ \hline
367  & 31  & 593  & 167 &  331 &  307 &  277  & 577  & 743   \\ \hline
311 &  67  & 191  & 151  & 281 &  47 &  101  & 619  & 439   \\ \hline
389  & 761  & 613 &  229 &  173 &  607 &  13  & 43  & 271   \\ \hline
421 &  563 &  241  & 557  & 317 &  337 &  673  & 751  & 113   \\ \hline
73  & 71  & 127  & 137 &  163 &  193 &  661 &  23  & 181   \\ \hline
409 &  571 &  691 &  61 &  83  & 251  & 179  & 233  & 877   \\ \hline
467 &  53  & 227 &  59 &  89 &  373 &  401  & 37  & 149   \\ \hline
19  & 547  & 809 &  521 &  131 &  41  & 659  & 503  & 491   \\ \hline
\end{tabular}
\caption{Prime pyramid of order 9 with weight 27325.}
\label{tab:pyramid9}
\end{table}

\section{Prime Cylinders}
\label{sec:cylinders}

\begin{defn}
A \emph{prime cylinder} of order $n$ with $k$ layers is a $n \times k$ matrix P of odd primes, such that for every $c$ and
$r>0$: $P(r,c)=P(r-1,c-1)+P(r-1,c)+P(r-1,c+1)$.
\end{defn}

Note that the columns wrap around and hence the term `cylinder'.
For example here is a prime cylinder of order 4 and 6 layers - the best found so far:

\begin{table}[!htpb]
\centering
\begin{tabular}{c}
$1091,3001,271,257$ \\
$4349,4363,3529,1619$ \\
$10331,12241,9511,9497$ \\
$32069,32083,31249,29339$ \\
$93491,95401,92671,92657$ \\
$281549,281563,280729,278819$
\end{tabular}
\end{table}

Since all the values below the first layer can be generated from previous values, a prime cylinder can be described using its
first layer only. So the above prime cylinder would be described as $(1091,3001,271,257)$. The \emph{weight} of a prime
cylinder is the sum of values in its first layer. When multiple prime cylinders have the same order and number of layers,
then we prefer the one with the smaller weight. Prime cylinders were originally introduced in \cite{Puzzle831}, but were
limited to $n=4$. Here we investigate other values of $n$. It turns out that prime cylinders of odd orders cannot have
more than two layers, so we focus on prime cylinders of even orders.
Table~\ref{tab:cyl} shows the best prime cylinders found for $n \leq 12$.


\begin{table}[!htpb]
\centering
\begin{tabular}{|c|c|c|c|}
\hline
$n$ & Prime Cylinder & Layers & Weight \\ \hline
4   & $(1091,3001,257,271)$   & 6      & 4620     \\ \hline       
6   & $(163,1109,307,1163,109,1307)$ & 6   &  4158    \\ \hline   
8   & $(67,541,23,137,109,193,389,431)$ & 5 & 1890     \\ \hline
10  & $(19,17,7,107,43,23,13,71,79,101)$ & 4 & 480  \\ \hline
12  & $(11,29,31,79,53,5,109,43,47,41,61,139)$ & 4 & 648  \\ \hline
\end{tabular}
\caption{Best prime cylinders found for $n \leq 12$.}
\label{tab:cyl}
\end{table}

\section{Goldbach Squares}
\label{sec:goldbach}

The famous Goldbach conjecture states that

\emph{Every even integer greater than 2 can be expressed as the sum of two primes}.

Although the conjecture has been verified up to $4 \times 10^{18}$~\cite{Goldbach}, a proof still remains elusive.
Here we investigate a problem related to the Goldbach conjecture:
can we place primes into a square such that every even number is generated as the sum of two adjacent cells?
This puzzle has been explored in \cite{Puzzle835}. More formally we have:

\begin{defn}
A \emph{Goldbach square} of order $n$ is a $n \times n$ matrix of odd primes (not necessarily unique) such that the sum
of any two adjacent cells is one of the even numbers from $6$ to $4+4n(n-1)$ inclusive and every even number
in this range appears exactly once.
\end{defn}

For example, here is a Goldbach square of order 3:

\begin{table}[!htpb]
\centering
\begin{tabular}{|c|c|c|}
\hline
7 & 5 & 3 \\ \hline
17 & 11 & 3 \\ \hline 
3 & 7 & 19 \\ \hline
\end{tabular}
\label{tab:Goldbach_example}
\end{table}
The sums across rows are:

\begin{equation}
\begin{split}
7+5 = 12,~~~5+3 = 8,\\
17+11 = 28,~~~11+3 = 14,\\
3+7 = 10,~~~7+19 = 26.
\end{split}
\end{equation}
 
The sums down columns are:

\begin{equation}
\begin{split}
7+17 = 24,~~~17+3 = 20,\\
5+11 = 16,~~~11+7 = 18,\\
3+3 = 6,~~~3+19 = 22.
\end{split}
\end{equation}

Notice that every even number from 6 to 28 appears exactly once. If there are multiple Goldbach squares for a given
$n$ then we prefer the one with the smallest sum of cells (\emph{weight}). Tables~\ref{tab:Goldbach} and~\ref{tab:Goldbach2}
show the best Goldbach squares that we found for $n \leq 10$.

\begin{table}[!htpb]
\centering
\begin{tabular}{|c|c|c|}
\hline
n & Goldbach Square & Weight \\ \hline

2 &
\begin{tabular}{|c|c|}
\hline
5 & 7 \\ \hline
3 & 3 \\ \hline
\end{tabular}
& 18 \\[4ex] \hline

3 &
\begin{tabular}{|c|c|c|}
\hline
3 & 5 & 7 \\ \hline
7 & 11 & 17 \\ \hline
19 & 3 & 3 \\ \hline
\end{tabular}
& 75 \\[6ex] \hline

4 &
\begin{tabular}{|c|c|c|c|}
\hline
5  & 11 &  13  & 5     \\ \hline
17 & 19  & 31 &  7     \\ \hline
3  & 29  & 11 &  3     \\ \hline
5  & 23  & 23 &  3     \\ \hline
\end{tabular}
& 208 \\[8ex] \hline

5 &
\begin{tabular}{|c|c|c|c|c|}
\hline
5 & 7 & 11 & 17 & 5 \\ \hline
31 & 31 & 23 & 53 & 3 \\ \hline
11 & 37 & 43 & 29 & 23 \\ \hline
3 & 47 & 31 & 29 & 17 \\ \hline
3 & 17 & 13 & 3 & 7 \\ \hline
\end{tabular}
& 499 \\[10ex] \hline

6 &
\begin{tabular}{|c|c|c|c|c|c|}
\hline
7 & 79 & 13 & 17 & 11 & 3 \\ \hline
41 & 17 & 71 & 43 & 31 & 19 \\ \hline
37 & 29 & 23 & 59 & 59 & 13 \\ \hline
61 & 7 & 47 & 53 & 3 & 31 \\ \hline
19 & 103 & 17 & 23 & 101 & 7 \\ \hline
5 & 13 & 3 & 3 & 5 & 5 \\ \hline
\end{tabular}
& 1078 \\[12ex] \hline

7 &
\begin{tabular}{|c|c|c|c|c|c|c|}
\hline
37 & 11 & 19 & 7 & 83 & 47 & 53 \\ \hline
5 & 17 & 109 & 61 & 19 & 17 & 53 \\ \hline
47 & 67 & 37 & 61 & 73 & 59 & 3 \\ \hline
11 & 7 & 71 & 11 & 83 & 89 & 7 \\ \hline
13 & 151 & 17 & 149 & 3 & 29 & 109 \\ \hline
3 & 3 & 107 & 13 & 137 & 5 & 3 \\ \hline
17 & 37 & 29 & 31 & 7 & 7 & 43 \\ \hline
\end{tabular}
& 2077 \\[15ex] \hline

\end{tabular}
\caption{Goldbach squares for $2 \leq n \leq 7$.}
\label{tab:Goldbach}
\end{table}

\begin{table}[!htpb]
\centering
\begin{tabular}{|c|c|c|}
\hline
n & Goldbach Square & Weight \\ \hline

8 &
\begin{tabular}{|c|c|c|c|c|c|c|c|}
\hline
17 & 23 & 157 & 13 & 43 & 37 & 31 & 67 \\ \hline
17 & 7 & 47 & 71 & 79 & 79 & 97 & 127 \\ \hline
3 & 71 & 131 & 31 & 61 & 127 & 37 & 73 \\ \hline
23 & 5 & 7 & 11 & 149 & 19 & 89 & 139 \\ \hline
41 & 53 & 113 & 41 & 3 & 3 & 127 & 17 \\ \hline
149 & 59 & 11 & 131 & 83 & 13 & 47 & 179 \\ \hline
71 & 29 & 3 & 5 & 31 & 59 & 3 & 7 \\ \hline
11 & 37 & 181 & 41 & 151 & 47 & 101 & 31 \\ \hline
\end{tabular}
& 3766 \\[17ex] \hline

9 &
\begin{tabular}{|c|c|c|c|c|c|c|c|c|}
\hline
101 & 47 & 13 & 107 & 97 & 131 & 113 & 59 & 149 \\ \hline
59 & 191 & 43 & 107 & 31 & 101 & 163 & 13 & 73 \\ \hline
107 & 29 & 19 & 3 & 7 & 173 & 109 & 139 & 7 \\ \hline
139 & 79 & 13 & 13 & 11 & 89 & 127 & 151 & 137 \\ \hline
23 & 173 & 37 & 31 & 83 & 7 & 7 & 13 & 61 \\ \hline
89 & 113 & 29 & 227 & 41 & 5 & 23 & 41 & 43 \\ \hline
101 & 179 & 5 & 3 & 127 & 73 & 29 & 41 & 199 \\ \hline
73 & 3 & 37 & 223 & 31 & 67 & 89 & 17 & 71 \\ \hline
113 & 3 & 151 & 61 & 163 & 103 & 103 & 19 & 107 \\ \hline
\end{tabular}
& 6187 \\[19ex] \hline

10 &
\begin{tabular}{|c|c|c|c|c|c|c|c|c|c|}
\hline
89 & 11 & 101 & 71 & 271 & 13 & 5 & 89 & 59 & 197 \\ \hline
107 & 227 & 23 & 251 & 67 & 19 & 3 & 211 & 29 & 47 \\ \hline
13 & 83 & 47 & 79 & 43 & 313 & 3 & 139 & 127 & 131 \\ \hline
3 & 223 & 11 & 163 & 7 & 13 & 181 & 181 & 37 & 61 \\ \hline
89 & 79 & 73 & 151 & 3 & 173 & 73 & 131 & 227 & 97 \\ \hline
199 & 149 & 179 & 97 & 41 & 7 & 127 & 233 & 71 & 47 \\ \hline
97 & 17 & 11 & 31 & 5 & 7 & 109 & 53 & 37 & 31 \\ \hline
257 & 89 & 251 & 101 & 131 & 19 & 37 & 241 & 67 & 7 \\ \hline
79 & 109 & 31 & 191 & 29 & 11 & 43 & 19 & 163 & 53 \\ \hline
127 & 163 & 181 & 89 & 179 & 23 & 59 & 5 & 47 & 19 \\ \hline
\end{tabular}
& 9212 \\[21ex] \hline

\end{tabular}
\caption{Goldbach squares for $8 \leq n \leq 10$.}
\label{tab:Goldbach2}
\end{table}

\section{Prime Matrices}
\label{sec:matrices}

\begin{defn}
A \emph{prime matrix} of order $n$ is a $n \times n$ matrix $P$ of odd primes, such that the sum of every odd number of elements in any 
straight line is prime. More formally, we have

\begin{equation}
\begin{split}
\sum_{0 \leq k \leq 2L} P(r+kd_r,c+kd_c) & ~~~\mbox{is prime for} \\
& \forall~r,c~\mbox{such that}~0 \leq r,c < n~\mbox{and} \\
& \forall~d_r,d_c~\mbox{such that}~(d_r,d_c) \in \{(0,1),(1,0),(1,1)\}~\mbox{and} \\
& \forall~L~\mbox{such that}~0 \leq 2L < n~\mbox{and} \\
& ~r+2Ld_r<n~\mbox{and}~c+2Ld_c<n.
\end{split}
\end{equation}
\end{defn}

We were able to find prime matrices up to order 7. For $n \leq 4$ we found optimal (smallest weight) prime matrices.
The results can be seen in Table~\ref{tab:matrix}. The lower bound on the optimal weight is the sum of the first $n^2$ odd primes.

\begin{table}[!htpb]
\centering
\begin{tabular}{|c|c|c|c|}
\hline
n & Prime Matrix & Weight & Lower Bound \\ \hline

3 &
\begin{tabular}{|c|c|c|}
\hline
5 & 19 & 13  \\ \hline  
3 & 17 & 23  \\ \hline 
29 & 11 & 7  \\ \hline 
\end{tabular}
& \multicolumn{2}{c|}{127} \\[7ex] \hline   

4 &
\begin{tabular}{|c|c|c|c|}
\hline
19 & 53 & 17 & 3    \\ \hline
47 & 23 & 31 & 29    \\ \hline
43 & 7 & 11 & 5    \\ \hline
41 & 13 & 59 & 37    \\ \hline
\end{tabular}
& \multicolumn{2}{c|}{438}  \\[8ex] \hline

5 &
\begin{tabular}{|c|c|c|c|c|}
\hline
101 & 107 & 61 & 109 & 41    \\ \hline
127 & 43 & 11 & 29 & 31    \\ \hline
5 & 17 & 37 & 13 & 59    \\ \hline
79 & 97 & 23 & 131 & 19    \\ \hline
47 & 53 & 7 & 67 & 89    \\ \hline
\end{tabular}
& 1403 & 1159 \\[10ex] \hline

6 &
\begin{tabular}{|c|c|c|c|c|c|}
\hline
73 & 131 & 113 & 109 & 61 & 227    \\ \hline
149 & 541 & 229 & 41 & 11 & 211    \\ \hline
379 & 491 & 419 & 349 & 139 & 53    \\ \hline
89 & 97 & 13 & 71 & 83 & 277    \\ \hline
193 & 59 & 137 & 307 & 127 & 29    \\ \hline
167 & 43 & 31 & 5 & 101 & 241    \\ \hline
\end{tabular}
& 5796 & 2582 \\[12ex] \hline

7 &
\begin{tabular}{|c|c|c|c|c|c|c|}
\hline
547 & 719 & 1117 & 983 & 1201 & 29 & 397\\ \hline 
691 & 827 & 103 & 1307 & 373 & 131 & 37\\ \hline 
53 & 457 & 503 & 7 & 419 & 73 & 557 \\ \hline
347 & 463 & 683 & 307 & 647 & 337 & 1433 \\ \hline
163 & 167 & 313 & 1013 & 127 & 1217 & 643 \\ \hline
367 & 677 & 787 & 107 & 193 & 653 & 13 \\ \hline
1223 & 1087 & 23 & 769 & 227 & 487 & 887 \\ \hline
\end{tabular} 
& 25891 & 5115 \\[15ex] \hline
\end{tabular}
\caption{Prime matrices for $3 \leq n \leq 7$.}
\label{tab:matrix}
\end{table}

\clearpage

\section{Conclusion and Future Work}

We have investigated a number of constructions that generate primes via the sum of primes. 
Some constructions are more efficient than others at generating primes. We can define
a construction's \emph{efficiency} as the number of primes it generates divided by the
number of primes used to construct the construction. Table~\ref{tab:eff} shows the greatest efficiency
achieved by each construction sorted from highest to lowest:

\begin{table}[!htpb]
\centering
\begin{tabular}{|c|c|c|}
\hline
\textbf{Construction} & \textbf{Order $n$} & \textbf{Efficiency} \\ \hline
Cyclic Prime Vector & 14 & 6 \\ \hline
Prime Vector & 23 & 5.26 \\ \hline
Prime Cylinder & 4 and 6 & 5 \\ \hline
Prime Matrix & 7 & 4 \\ \hline
Prime Stair & 15 & 3.27 \\ \hline
Prime Tuple & 7 & 2.57 \\ \hline
Prime Pyramid & 9 & 1.04 \\ \hline
\end{tabular}
\label{tab:eff}
\end{table}

Many questions remain unresolved:

\begin{itemize}
\item What are the optimal prime vectors for $15 \leq n \leq 23$ ?
\item Is there a prime vector of order 24 ?
\item What are the optimal cyclic prime vectors for $11 \leq n \leq 14$ ?
\item Is there a cyclic prime vector of order 15 ?
\item What are the optimal prime tuples for $n=7, 13, 19$ ?
\item What are the optimal prime stairs for $12 \leq n \leq 15$ ?
\item Is there a prime stair of order 16 ?
\item What is the optimal prime pyramid of order 9 ?
\item Is there a prime pyramid of order 10 ?
\item What are the optimal prime cylinders for $n \leq 12$ ?
\item Is there a Goldbach square of order 11 ?
\item What are the optimal prime matrices for $5 \leq n \leq 7$ ?
\item Is there a prime matrix of order 8 ?

\end{itemize}

\section{Acknowledgements}
We are very grateful to Trevor Tao for reviewing the paper and providing great suggestions.
We also want to thank Carlos Rivera for running primepuzzles.net that inspired many of the constructions in this paper.

\vspace{1cm}
\bibliographystyle{plain}
\bibliography{shortbib}

\begin{thebibliography}{1}

\bibitem{Progression}
{Primes in arithmetic progression}.
\newblock \url{https://en.wikipedia.org/wiki/Primes_in_arithmetic_progression}.

\bibitem{Progression2}
Jens~Kruse Andersen.
\newblock {Primes in arithmetic progression records}.
\newblock \url{http://primerecords.dk/aprecords.htm}.

\bibitem{Goldbach}
Tom\'{a}s~Oliveira e~Silva.
\newblock {Goldbach conjecture verification}.
\newblock \url{http://www.ieeta.pt/~tos/goldbach.html}.

\bibitem{Magic}
Harvey Heinz.
\newblock {Prime Magic Squares}.
\newblock \url{http://recmath.org/Magic%20Squares/primesqr.htm}.

\bibitem{Puzzle831}
Carlos Rivera.
\newblock {Puzzle 831. Rings of primes}.
\newblock \url{http://primepuzzles.net/puzzles/puzz_831.htm}.

\bibitem{Puzzle835}
Carlos Rivera.
\newblock {Puzzle 835. Goldbach squares}.
\newblock \url{http://primepuzzles.net/puzzles/puzz_835.htm}.

\bibitem{OEIS}
Neil Sloane.
\newblock {The On-Line Encyclopedia of Integer Sequences}.
\newblock \url{http://oeis.org}.

\bibitem{Array}
Eric~W. Weisstein.
\newblock {Prime Array}.
\newblock \url{http://mathworld.wolfram.com/PrimeArray.html}.

\bibitem{Magic2}
Eric~W. Weisstein.
\newblock {Prime Magic Square}.
\newblock \url{http://mathworld.wolfram.com/PrimeMagicSquare.html}.

\end{thebibliography}

\end{document}